\documentclass{article}

\usepackage{amsmath, amssymb, amsthm, amscd}
\usepackage{graphics}

\newtheorem{lemma}        {Lemma}
\newtheorem{lowerbound}   {Lower bound}
\newtheorem{proposition}  {Proposition}

\theoremstyle{definition}

\newtheorem{definition}   {Definition}

\newtheorem{problem}      {Problem}

\newcommand{\slashpoly}    {\ensuremath{_\mathrm{/poly}}}
\newcommand{\ppoly}        {\ensuremath{\mathcal P\slashpoly}}

\title{Lower bounds for some decision problems over $\mathbb C$}
\author{Gregorio Malajovich\footnote{
Departamento de Matem\'atica Aplicada, Universidade Federal
do Rio de Janeiro. Caixa Postal 68530, CEP 21945, Rio de Janeiro,
RJ, Brasil. e-mail: gregorio@labma.ufrj.br.
On leave at MSRI, 1000 Centennial Drive, Berkeley CA
94720-5070. e-mail: gregorio@msri.org
}}
\date{January 22, 1999}
\begin{document}
\maketitle
\begin{abstract}
    Lower bounds for some explicit decision problems over the
    complex numbers are given.
\end{abstract}

\section{Introduction}
	This paper is about lower bounds for certain decision
	problems over $\mathbb C$. (See ~\cite{BCSS98} for the
	model of computation and for background). In particular,
	we will provide lower bounds for the complexity of
	deciding, given $x$, if $p^d(x)=0$ for some explicit
	polynomials $p^d$.
\par
	A related problem is to give lower bounds for the 
	evaluation of explicit polynomials. This
	has been an active subject of research since
	~\cite{STRASSEN74}.
	See~\cite{BCS} for modern developments and for bibliographical remarks.
	More recent results appeared in~\cite{AHMMP} and~\cite{AM98}.
\par
	Most of those bounds use the Ostrowsky model of computation
	~(\cite{BCS} page 6): 
	sum and multiplication by an algebraic constant are free,
	and the complexity of a computation for polynomial $f(x)$
	is the number of non-scalar multiplications, i.e., of 
	multiplications of two polynomials in the variable $x$.
	For instance, Horner rule for a degree $d$ polynomial 
	requires $d$ non-scalar multiplications.
\par
	All those bounds apply trivially to the complexity of
	evaluating polynomials by a `machine over $\mathbb C$'
	as defined in~\cite{BCSS98}, or to the 
	(multiplicative-branching) complexity of a 
	computation tree for evaluating the same polynomial.
\par
	Little is known, however, about the application of those
	bounds to decision problems (Over $\mathbb C$, in the sense
	of ~\cite{BCSS98}, or by a decision tree as in~\cite{BCS},
	Definition (4.19) page 115. In this definition, each node
	of a computation tree can perform one algebraic operation
	or comparison, and therefore a natural measure of complexity
	is the depth of the tree).
\par
	In this paper, only decision problems of the form below will
	be considered: let $X \subseteq \mathbb N \times \mathbb C$,
	and let $X_d = \{ x \in \mathbb C : (d,x) \in X \}$. 
	Typically, $d$ is the problem size and $\# X_d \le d$.
	One can think of $X$ as the disjoint union of the zero-set
	of a family of polynomials of degree $\le d$, where $d \in
	\mathbb N$. The two following forms of a decision problem
	are natural in this setting:
\begin{problem} \label{pnu}
	For any fixed $d$, decide wether $x \in X_d$.
\end{problem}

\begin{problem} \label{pu}
	Decide wether $(d,x) \in X$
\end{problem}
\medskip
\par
	Problem~\ref{pnu} is non-uniform, in the sense that we allow
	for a different machine over $\mathbb C$ or a different decision
	tree to be used for each value of $d$. However, we want a 
	bound on the running time or on the multiplicative complexity
	of the tree, as a function of $d$. 
\par
	Problem~\ref{pu} is uniform. It is harder than Problem~\ref{pnu},
	in the sense that it cannot be solved by a decision tree, since
	$\# X_d$ can be arbitrarily large. It requires
	a machine over $\mathbb C$, that will eventually branch according
	to the value of $d$.
\par
	Lower bounds for Problem~\ref{pnu} are also lower bounds
	for Problem~\ref{pu}. 
\par
	A trivial, topological lower bound for Problems~\ref{pnu} 
	and~\ref{pu} when $\# X_d = d$ is $\log_2 d$. Sharper known
	bounds come from the `Canonical Path' argument, see
	~\cite{BCSS98} section 2.5:
	Let $f$ be a univariate polynomial. The complexity of deciding
	$f(x)=0$ is bounded below by the minimum of the complexity
	of evaluating $g(x)$, where $g$ ranges over the non-zero
	multiples of $f$.
\par
	If one
	assumes some property of $f$ that propagates to
	its multiple $g$, then one eventually obtains a
	sharper, non-trivial lower bounds.
\par
	In Lemma ~\ref{lower1} below, we will give
	conditions on the roots of $f$ that will provide lower bounds 
	for the evaluation of $g$. Essentially, we will require a subset
	of the roots to be rapidly growing. This will imply a rapid
	growth property for the coefficients of $g$. Then, the results
	of~\cite{AHMMP, AM98}
	imply a lower bound for the complexity of
	evaluating $g$. Thus 
	we will be able to construct specific
	polynomials that are hard to decide in the {\em non-uniform}
	sense, viz.

\begin{lowerbound} \label{low1}
   The set $X = \{ (d,x) \in \mathbb Z \times \mathbb C :
   x = 2^{2^{di}}, 0 \le i \le d$, cannot
   be solved in time polylog($d$) in the setting of Problem~\ref{pnu}.
\end{lowerbound}

\begin{lowerbound} \label{low2}
   The set $Y = \{ (d,x)  \in \mathbb Z \times \mathbb C :
   p^{d}(x)=0 \}$, where $p^d(t) = \sum _{i=0}^d 2^{2^{d (d-i)}} t^i$,
   cannot be decided in time polylog($d$) in the setting of Problem~\ref{pnu}.
\end{lowerbound}
\par
	In a more classical computer-science language, we can define the
	input size of some $(d,x)$ as $\log d$. This means that the integer
	$d$ is represented in binary notation, while variable $x$ can
	contain an arbitrary complex number. In that case, `time 
	polylog($d$) in the setting of Problem~\ref{pnu}' can be
	refrased as $\ppoly$. The lower bounds above become now: $X
	\not \in \ppoly$ and $Y \not \in \ppoly$.
\par
	Non-uniform lower bounds ~\ref{low1} and ~\ref{low2} can be
	compared to the following easier, uniform lower bound:

\begin{lowerbound}\label{low3}
   The set $Z = \{ (d,x)  \in \mathbb Z \times \mathbb C :
   q^d(x)=0 \}$, where $q^d(t) = \sum _{i=0}^d 2^{2^i} t^i$,
   cannot be decided in time polylog($d$) in the setting of Problem~\ref{pu}.
\end{lowerbound}
\par
	This means that the set $Z$, where $d$ is represented in binary
	notation and $x$ is a complex number, does not belong to
	$\mathcal P$ over $\mathbb C$.
\medskip
\par
	This work was written while the author was visiting the
	Mathematical Sciences Research Institute, 
	Berkeley, CA. Thanks to Pascal Koiran, Jos{\'e} Luis
	Monta\~na, Luis Pardo and Steve Smale for their 
	suggestions and comments.

\section{Background and notations}
\begin{definition}
	Let $K \subset L$ be finite algebraic extensions of $\mathbb Q$.
	Let $\nu$ be a valuation in $M_K$. Then we extend the notation $\nu$
	to $L$ by:
\[
	\nu(x) = \frac{ \sum_{\mu} n_{\mu} \mu(x) }{\deg[L:K]}
\]
	where the sum ranges over all the valuations $\mu$ of $L$ that
	are `above' $\nu$, and where $n_{\mu}$ is the `local degree' of
	$L:K$. The local degree is defined as 
	$n_{\mu} = \deg[ L_{\mu} : K_{\nu} ]$, where $K_{\nu}$ is the 
	completion of $K$ under the metric induced by the absolute
	value $|.|_{\nu}$. 
\end{definition}

	Recall that for $x \in K$, $\deg[L:K] 
	\nu(x) = \sum_{\mu} n_{\mu} \nu(x)$.
	The case $K = \mathbb Q$ is an immediate consequence of
	Corollary 2 of Theorem 1 in Chapter II, p. 39 of \cite{LANG86}.

\begin{definition}
	Let $g$ be a polynomial with algebraic coefficients in some
	extension $K$ of $\mathbb Q$. Let $\nu$ be a valuation in $M_K$.
	The {\em Newton diagram} of $g$ at $\nu$ is the (lower) convex
	hull of the set $\{ (i, \nu(g_i)), i=0 \cdots d \}$.
\end{definition}

	The basic property of Newton diagrams used here is the following.
\begin{proposition}\label{prop1}
	Suppose that $\zeta_1, \cdots, \zeta_d$ are the roots of a univariate
	polynomial $g \in K[x]$.
	Let the roots of $g$ be ordered so that
\[
	\nu(\zeta_1) \ge \cdots \ge \nu(\zeta_d)
\]
	and let the increasing sequence $i_j$ assume the values $0$,
	$d$ and all the values of $i$ where:
\[
	\nu(\zeta_i) > \nu(\zeta_{i+1})
\]
\par
	Then the sharp corners of the Newton diagram are precisely the
	points of the form $(i_j, \nu(g_{i_j}))$ for 
	all $j$.
\par
	Moreover, the slope of the segment $[(i_{j-1}, \nu(g_{i_{j-1}})), 
	(i_{j}, \nu(g_{i_{j}}))]$ is precisely $- \nu(\zeta_{i_{j}})$.
\end{proposition}
\begin{proof}[Proof of Proposition~\ref{prop1}]
	The proof uses the following property of valuations: 
	$\nu (\sum x_i) \ge \min \nu (x_i)$. Furthermore, when that 
	minimum is attained in only one $x_i$, we have equality.
\par
	Let $i_{j-1} < k < i_j$. Writing
\begin{eqnarray*}
	g_{i_{j-1}} &=& g_d \sigma_{d-i_{j-1}}(\zeta_1, \cdots, \zeta_d) \\
	g_{k} &=& g_d \sigma_{d-k}(\zeta_1, \cdots, \zeta_d) \\
	g_{i_{j}} &=& g_d \sigma_{d-i_{j}}(\zeta_1, \cdots, \zeta_d) \\
\end{eqnarray*}
\par
	one can pass to the valuation by:
\begin{eqnarray*}
	\nu(g_{i_{j-1}}) &=& \nu(g_d) + \nu(\zeta_{i_{j-1}+1}) 
	+ \cdots + \nu(\zeta_d) \\
	\nu(g_{k})       &\ge& \nu(g_d) + \nu(\zeta_{k+1}) 
	+ \cdots + \nu(\zeta_d) \\
	\nu(g_{i_{j}})   &=& \nu(g_d) + \nu(\zeta_{i_{j}+1}) 
	+ \cdots + \nu(\zeta_d) \\
\end{eqnarray*}
\par
	Subtracting, one obtains:
\begin{eqnarray*}
	\nu(g_{i_j})       
	- \nu(g_{i_{j-1}}) 
	&=& -\nu(\zeta_{i_{j-1}+1}) - \cdots - \nu(\zeta_{i_{j}}) \\
	&=& - (i_j - i_{j-1}) \nu(\zeta_{i_j}) \\
	\nu(g_{i_j})       
	- \nu(g_{k}) 
	&\le& -\nu(\zeta_{k+1}) - \cdots - \nu(\zeta_{i_{j}}) \\
	&\le& - (i_j - k) \nu(\zeta_{i_j}) \\
\end{eqnarray*}
\par
	This concludes the proof.	
\end{proof}
\medskip
\par
\section{Uniform lower bounds}

	We can now prove Lower Bound~\ref{low3}.
\begin{proof}[Proof of Lower bound~\ref{low3}]
\centerline{\resizebox{10cm}{6cm}{\includegraphics{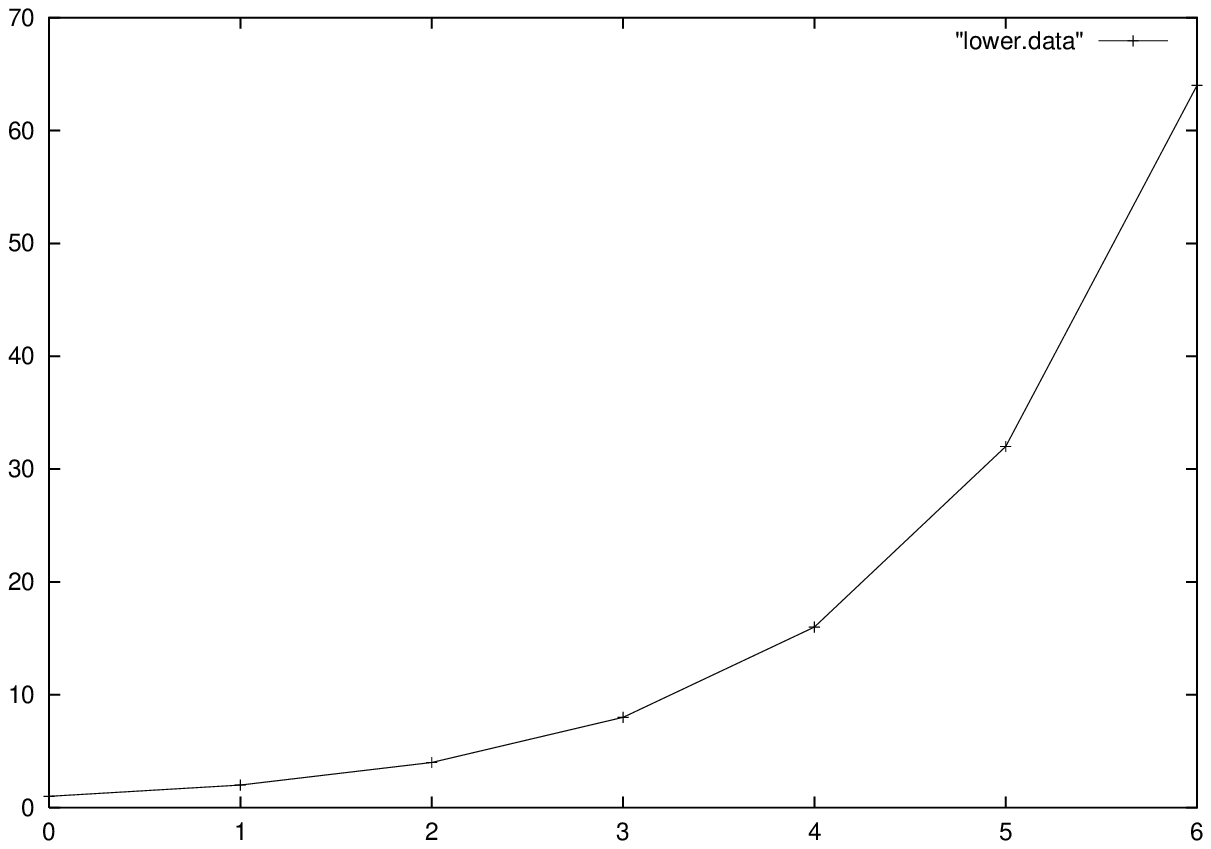}}}

	The Newton diagram of $q^{d}$ at 2 is $\{ (i, 2^i): 0 \le i \le d\}$.
	(This latest set is convex, since the points lie on the curve $y=2^x$
	and this curve is convex). Therefore, there is a unique root $\zeta$
	of $q^{d}$ that minimizes $\nu(\zeta)$. 
\par
	Since $q^{d}_{d-1} = (-\sum \zeta_i) q^{d}_d$, where the sum
	ranges over all the roots, we have:
\[
	\nu_2 (q^{d}_{d-1}) =  \nu_2 (q^{d}_{d}) + \min \nu_2(\zeta_i)
	= \nu_2 (q^{d}_{d}) + \nu_2(\zeta)
\]
	Replacing by the actual values of the coefficients, one gets:
\begin{equation}\label{eq1}
	\nu_2(\zeta) = -2^{d-1} 
\end{equation}
\medskip
\par
	Now, suppose that there is a machine $M$ that decides $q^{d}(t) = 0$
in time polylog($d$). One can assume without loss of generality that this
machine has no constant but $0$ and $1$. Let its running time be bounded by
$T=a (\log d)^b$.
\par
	Let us fix $d > 2+T^2$. We will derive a contradiction.
\par
	Let $g$ be the polynomial defining the canonical path (recall that
	$d$ is fixed now, so this is the path followed by generic $t \in
	\mathbb C$). It can be
	computed in time $\le T^2$, so we have the following bounds:
\begin{eqnarray*}
	\deg g &\le 2^{T^2} \\
	0\le \nu_2 (g_p) &\le 2^{T^2}
\end{eqnarray*}
\par
	Since $\zeta$ is also a root of $g$, there are coefficients
	$g_i$ and $g_j$, $i \ne j$, such that:
\begin{equation}\label{eq2}
	(j-i) \nu_2(\zeta) = \nu_2(g_i) - \nu_2 (g_j)
\end{equation}
	\par
	Thus, $|\nu_2(\zeta)| \le |\nu_2(g_i)|+|\nu_2(g_i)|$.
	This implies:
\[
	|\nu_2(\zeta)| \le {2^{1+T^2}} < 2^{d-1}
\]
	Replacing by equation~\ref{eq1}, one obtains $2^{d-1}<2^{d-1}$,
	a contradiction.
\end{proof}

\section{Non-uniform lower bounds}

\begin{lemma}\label{lower1}
	Let $g=g(t)$ be a degree $D$ polynomial with algebraic coefficients.
Let $\nu$ be a (non-archimedian) valuation of $K=\mathbb Q[g_0, \cdots, g_D]$.
	Let $\xi_1, \cdots \xi_D$ be the roots of $g$, and assume they
	are ordered in such way that:
\[
	\nu(\xi_1) \ge \cdots \ge \nu(\xi_D)
\]
	Suppose that there is a subsequence $\zeta_j = \xi_{i_j+1}$,
	$j=1 \cdots d$, such that the following holds:
	\begin{enumerate}
	\item $\nu(\zeta_d) \ge 1$ 
	\item $\nu(\zeta_{j}) \ge 2 (i_{j+1} - i_{j}) 
	       \ \nu(\zeta_{j+1})$, for $0 \le j \le d-1$.
	\end{enumerate}
\medskip
\par
	Then $g$ cannot be evaluated in less than 
\[
	L \ge \sqrt{ \frac{d}{28 \log_2 D+1} }
\]
	multiplications.
\end{lemma}

\begin{proof}[Proof of Lemma~\ref{lower1}]
	We can assume without loss of generality that the ordering
	of the $\xi_i$ satisfies:
\[
	\cdots \xi_{i_j} < \xi_{i_j+1} = \zeta_j \le \xi_{i_j +2} \cdots
\]
	For $j \in \{1, \cdots , d-1 \}$ we have:
\[
	\nu ( g_{i_j} ) - \nu (g_{i_{j+1}} ) = 
	\nu(\xi_{i_j+1}) + \cdots + \nu(\xi_{i_{j+1}})
\]
	Hence, using $\nu(\xi_{i_{j+1}}) > \nu(\zeta_d) \ge 1$:
\[
	\nu(\zeta_j) 
	\le \nu ( g_{i_j} ) - \nu (g_{i_{j+1}} ) \le 
	(i_{j+1} - i_{j}) \nu (\zeta_j) 
\]
	By the same argument, for $j \in \{0, \cdots , d-2\}$:
\[
	\nu(\zeta_{j+1}) 
	\le \nu ( g_{i_{j+1}} ) - \nu (g_{i_{j+2}} ) \le 
	(i_{j+2} - i_{j+1}) \nu (\zeta_{j+1}) 
\]
\par
	Hence,
\[
	\frac{ \nu ( g_{i_j} ) - \nu (g_{i_{j+1}} ) }
	     { \nu ( g_{i_{j+1}} ) - \nu (g_{i_{j+2}} ) }
	\ge
	\frac {\nu(\zeta_{j})} {(i_{j+1} - i_{j}) \nu (\zeta_{j+1})}
	\ge 2
\]
\par
	Set $G_j = \nu(g_{i_j})$ for $j=0, \cdots, d-1$. We know that
	the $G_j$ are such that $|G_{j+1} - G_j| < 
	\frac{1}{2} |G_j - G_{j-1}|$.
	Hence 
\[
	\# \{ \sum s_j G_j, s_j \in \{0;1\} \} = 2^d
\]
\par
	Hence:
\[
	\# \{ \nu (\prod_{s \in S} g_{s}) , S \subset \{0, \cdots, D\} \} 
	      \ge 2^d
\]
\par
	and hence	
\[
	\mu(g) = \# \{ \sum_{S \subset \{0, \cdots, D\}}
	\theta_S \prod_{s \in S} g_{s} , \theta_S 
	\in \{0;1\} \} \ge 2^{2^d}
\]
\par
	By Lemma~1 in~\cite{AHMMP} or by Lemma 4 in~\cite{AM98}, 
\[
	\mu(g) \le 2^{ (D+1)^{28 L^2} }
\]
	and hence, taking logs:  
\[
	(D+1)^{28 L^2} \ge 2^d 
\]
	Taking logs again:
\[
	28 L^2 \ge \frac{d}{\log_2 D+1}
\]
	and hence:
\[
	L \ge \sqrt{ \frac{d}{28 \log_2 D+1} }
\]
\end{proof}
\par
	Note: Lemma~1 in~\cite{AHMMP} is slightly more general than
	Lemma~4 in~\cite{AM98}. However, using Lemma~4 in~\cite{AM98}
	it is possible to replace all the appearances of the number 28 in the
	statement and proof of Lemma~\ref{lower1} above by the number 21. 
\par
\begin{proof}[Proof of Lower Bound~\ref{low2}]
	We see from its Netwon diagram that the polynomial $p$ has
distinct roots $\zeta_1, \cdots \zeta_d$ with:
\[
	\nu_2 ( \zeta_i ) = 2^{d (d-i+1)} - 2^{d (d-i)}
	=
	2^{d (d-i)} (2^d - 1) 
\]
\par
	So we have $\nu_2 (\zeta_d) = 2^d - 1 > 1$, and
\begin{equation}\label{rationu}
	\nu_2(\zeta_i) / \nu_2(\zeta_{i+1}) = 2^d
\end{equation}
\par
	Assume that there are $a, b$ such that for each $d$, there
is a machine $M$ over $\mathbb C$ deciding $p(t)=0$ in time
$T= a (\log d)^b$. Its generic path is defined by a polynomial
$g(t)$ of degree $\le 2^T$.
\par
	Let us fix $d > 28 (T+1) T^2$. In particular $d \ge T+1$.
	We are in the conditions of Lemma~\ref{lower1},
where $D=2^T$. From that Lemma, it follows that
\[
	T \ge \sqrt{ \frac{d}{28 \log_2 2^T +1} }
	  \ge \sqrt{ \frac{d}{28 (T+1)} }
\]
	Hence,
\[
	28 T^2 (T+1) \ge d 
\]
	contradicting our choice of $d$.
\end{proof}
	Equation ~(\ref{rationu}) holds trivially in the proof of 
	Lower bound~\ref{low1}. The rest of the proof is verbatim
	the same.
\bibliographystyle{plain}
\bibliography{lower}

\end{document}